\documentstyle[12pt,a4]{article}

\def\rixrel#1#2{\mathrel{\mathop{#2}\limits_{#1}}}

\def\p{\sf p}
\def\q{\sf q}
\def\m{\sf m}

\begin{document}

\title{Towards a theory of arithmetic degrees}
\vspace{10mm}


\author{Chikashi Miyazaki\thanks{This author would like to thank
 Massey
University for
financial support and the Department of Mathematics for its friendly
atmosphere while
writing this paper.}\\Nagano National College of Technology\\716
 Tokuma,
Nagano 381\\
Japan\\e-mail: miyazaki@ei.nagano-nct.ac.jp\\\
\and Wolfgang Vogel\\
Department of Mathematics\\Massey University\\Palmerston North\\New
 Zealand
\\ e-mail: W.\@Vogel@massey.ac.nz}

\date{}
\maketitle

\vspace{3cm}

\section{Introduction}

The aim of this paper is to start a systematic investigation of the
arithmetic degree of projective
schemes as introduced in \cite{BM}. One main theme concerns itself
 with the
behaviour of
this arithmetic degree under hypersurface sections, see Theorem 2.1.

The classical intersection theory only considers the top-dimensional
 (or
isolated) primary
components. However, the notion of arithmetic degree involves the
 new concept of
length-multiplicity of embedded primary ideals as considered in \cite{BM},
\cite{EH},
\cite{H}, \cite{K}, \cite{STV}. Therefore it is much harder to control
 the
arithmetic
degree under a hypersurface section than in the case for the classical
degree theory.
We describe in \S 3 an upper bound for the arithmetic degree in terms of
the Castelnuovo-Mumford regularity, see Theorem 3.1. In addition, we
generalize Bezout's theorem via iterated hypersurface sections, see
Theorem 4.1. We conclude in \S 5 by studying two examples.

\section{Arithmetic degree and hypersurface sections}

Before stating our main result of this section we need to introduce
 the concept
of the length-multiplicity of (embedded) primary components.

Let $K$ be an arbitrary field and $S$ the polynomial ring
$K[x_0,\cdots,x_n]$.
Let ${\m}=(x_0,\cdots,x_n)$ be the homogeneous maximal ideal
 of $S$.
Let $I$ be a
homogeneous ideal of $S$.
\vspace{5mm}

\noindent{\bf Definition (\cite{BM}):}
Let $\p$ be a homogeneous prime ideal belonging to $I$. For
 a primary
decomposition
$I=\cap {\q}$ we take the primary ideal $\q$ with
$\sqrt{{\q}}={\p}$.
Let $J$ be the intersection of all primary components of $I$ with
associated prime ideals
${\p}_1$ such that
${\p}_1\rixrel{\neq}{\subset}{\p}$. If
the prime ideal
$\p$ is an isolated component of $I$, we set $J=S$. We define
 the
length-multiplicity
of $\q$ denoted by $\mbox{mult}_I({\p})$, to be the length
 $\ell$
of a maximal
strictly increasing chain of ideals
$${\q}\cap J =: J_\ell\subset J_{\ell-1}\subset\cdots\subset
J_1\subset J_0 := J$$
where $J_k$, $1\leq k\leq \ell-1$, equals ${\q}_k\cap J$ for
 some
$\p$-primary
ideal ${\q}_k$.

Despite the non-uniqueness of embedded components, the number
$\ell=\mbox{mult}_I({\p})$ is
well-defined.

\vspace{5mm}

For a finitely generated graded $S$-module $M$, let $H(M,\ell)$
 be the
Hilbert function of $M$ for all
integers $\ell$, that is, $H(M,\ell)$ is the dimension of the vector
 space
$[M]_\ell$
over $K$. It is well-known that the Hilbert function of $M$ is a
 polynomial
in $\ell$
for $\ell$ large enough. We denote this polynomial by $P(M,\ell)$.

We set $\Delta(H(M,\ell))=H(M,\ell)-H(M,\ell-1)$, $\Delta^0
H(M,\ell)=H(M,\ell)$,
 and $\Delta^r(H(M,\ell))=
\Delta^{r-1}(\Delta H(M,\ell))$ for all integers $r\geq 2$. Moreover,
 we set
$$\Delta_\tau(H(M,\ell))=H(M,\ell)-H(M,\ell-\tau)$$
for all integers $\tau\geq 1$. Further, the Hilbert polynomial $P(M,\ell)$
is written as
$$P(M,\ell) = \frac{e}{d!}\ell^d + (\mbox{lower order terms}),
\quad e\neq 0.$$
Then we define $h$-dim $M=d$ (homogeneous dimension) and degree
 of $M$ by
$\deg M := e$. Also, we write, for any ideal $I$ of $S$, $\dim
 I$ and $\deg I$
for $h$-dim $S/I$ and $\deg S/I$ respectively. In case $P(M,\ell)=0$,
 we
define $h$-dim $M=-1$ and
$\deg M=\sum_{\ell\in{\bf Z}} H(M,\ell)$. In particular, $\dim
 {\m}=-1$ and $\deg {\m}=1$.
\vspace{5mm}

\noindent{\bf Definition (\cite{BM}):}
For an integer $r\geq -1$, we define
\begin{eqnarray*}
\mbox{arith-deg}_r (I) & = &
\sum_{\stackrel{{\p} \mbox{\ \scriptsize is a prime
ideal}}{\mbox{\scriptsize such that}
\dim{\p}=r}} \mbox{mult}_I({\p})\cdot\deg {\p}\\
& = & \sum_{\stackrel{{\p}\in \mbox{\scriptsize Ass }
S/I}{\mbox{\scriptsize such that }
\dim{\p}=r}}
\mbox{mult}_I(\p)\cdot\deg{\p}
\end{eqnarray*}
\vspace{3mm}

\noindent{\bf Definition (\cite{H}):}
Let $I$ be a homogeneous ideal of $S$. Let $r$ be an integer with
 $r\geq -1$. We
define the ideal $I_{\geq r}$ as the intersection of all primary
 components
$\q$
of $I$ with $\dim{\q}\geq r$.
\vspace{5mm}

\noindent{\bf Remark :}
Although primary decomposition is not uniquely determined, the
ideal $I_{\ge r}$ does not depend on the choice of
primary decomposition of $I$ and is again a homogeneous ideal.
\vspace{5mm}

The aim of this section is to prove the following theorem.
\vspace{5mm}

\noindent{\bf Theorem 2.1:}{\em\ Let $r$ be an integer with $r\geq
 0$. Let
$I$ be
a homogeneous ideal of $S$. Let $F$ be a homogeneous polynomial of
 $S$ with
$\mbox{degree }(F)=\tau\geq 1$. Assume that $F$ does not belong
 to any
associated
prime ideal $\p$ of $I$ with $\dim {\p}\geq r$. Then we
 have
$$\mbox{arith-deg}_{r-1}(I,F) - \mbox{arith-deg}_{r-1}(I_{\geq
r+1},F)\geq\tau\cdot
\mbox{arith-deg}_r (I)$$
and the equality holds if and only if $F$ does not belong to any
 associated
prime
ideal $\p$ of $I$ with $\dim{\p}=r-1$.}
\vspace{5mm}

\noindent{\bf Corollary 2.2:}{\em\  Under the above condition,
$$\mbox{arith-deg}_{r-1}(I,F)\geq\tau\cdot\mbox{arith-deg}_r(I)$$
and the equality holds if and only if $F$ does not belong to any
 associated
prime
ideal $\p$ of $I$ with $\dim{\p}=r-1$ and the ideal $(I_{\geq
r+1},F)$ has
no associated prime ideals of dimension $(r-1)$.}
\vspace{5mm}

\noindent{\em Proof.} Corollary 2.2 follows immediately from Theorem
 2.1.
\vspace{5mm}

We note that Corollary 2.2 and Lemma 3 of \cite{KS} yield Theorem 2.3 of
 \cite{STV}.

\vspace{5mm}

We want to consider generic hyperplane sections.
The following useful lemma is obtained from \cite{BF}, (4.2)
and \cite{F}, (5.2), which was pointed out to us by H. Flenner.

\vspace{5mm}

\noindent{\bf Lemma 2.3:}{\em\
Let $I$ be a homogeneous ideal of $S$. We set $A=S/I$. Let $h=0$
 be the
defining
equation of a generic hyperplane of ${\bf P}_K^n$ ($K$: infinite
 field).
Then we have
$$\mbox{Ass}(A/h)\verb+\+\{{\m}\}\rixrel{=}{\subset}
\bigcup_{{\p}\in
\mbox{\scriptsize Ass } A}
\mbox{Min} (A/({\p},h)),$$
where $\mbox{Min}(A/({\p},h))$ is the set of minimal primes
 belonging to
$({\p},h)$.}
\vspace{5mm}

\noindent{\bf Corollary 2.4:}{\em\
Let $r$ be an integer $\geq 1$. Let $H$ be a generic hyperplane in
 ${\bf
P}_K^n$, given
by $h=0$. Then we have:
$$\mbox{arith-deg}_r (I) = \mbox{arith-deg}_{r-1} (I,h).$$
}
\vspace{5mm}

\noindent{\em Proof.} Corollary 2.4 follows immediately from Theorem
 2.1 and Lemma 2.3.
\vspace{5mm}

We note that Corollary 2.4 is stated in \cite{BM}, page 33 without
proof.

\vspace{5mm}
Before we turn to the proof of Theorem 2.1, two technical results
 are
needed. First we
state a more or less known result describing
 a
different
characterization of the arithmetic degree, which is purely algebraic
 and,
in fact,
serves as the definition in \cite{H}.
\vspace{5mm}

\noindent{\bf Lemma 2.5:}{\em\
Let $r$ be a non-negative integer. Let $I$ be a homogeneous ideal
 of $S$.
Then we have
$$\mbox{arith-deg}_r (I) = \Delta^r (P(S/I,\ell) - P(S/I_{\geq
 r+1},\ell))$$
for all integers $\ell$.}
\vspace{5mm}

\noindent{\bf Lemma 2.6:}{\em\
Let $r$ be an integer with $r\geq 1$. Let $I$ be a homogeneous ideal
 of $S$
and $F$
a homogeneous polynomial of $S$ with $\deg(F)=\tau\geq 1$. Assume
 that $F$
does not
belong to any associated prime ideal ${\p}$ of $I$ with $\dim
 {\p}\geq r+1$. Then we have
\begin{eqnarray*}
\lefteqn{\mbox{arith-deg}_{r-1}(I,F)-\Delta^{r-1}P(S/(I_{\geq r+1},F),\ell)
 +
\Delta^{r-1}P(S/(I,F)_{\geq r},\ell)} \\
& & \qquad\qquad= \tau\cdot\mbox{arith-deg}_r (I) + \Delta^{r-1}
P([0:F]_{S/I},\ell-\tau)
\end{eqnarray*}
for all integers $\ell$.}
\vspace{5mm}

\noindent{\em Proof.}
   From the exact sequences
$$0\rightarrow [0:F]_{S/I}(-\tau)\rightarrow
S/I(-\tau)\stackrel{F}{\rightarrow}S/I
\rightarrow S/(I,F)\rightarrow 0$$
and
$$0\rightarrow S/I_{\geq r+1}(-\tau)\stackrel{F}{\rightarrow}S/I_{\geq
r+1}\rightarrow
S/(I_{\geq r+1},F)\rightarrow 0,$$
we have
$$\Delta_\tau P(S/I,\ell)=P(S/(I,F),\ell)-P([0:F]_{S/I},\ell-\tau)$$
and
$$\Delta_\tau P(S/I_{\geq r+1},\ell)=P(S/(I_{\geq r+1},F),\ell)$$
for all integers $\ell$.
Note that
 $P(S/I,\ell)   -   P(S/I_{\ge r+1},\ell)$
is a numerical polynomial of degree $r$, see (2.5).
Thus we have
\begin{eqnarray*}
\lefteqn{\mbox{arith-deg}_{r-1}(I,F)-\Delta^{r-1}P(S/(I_{\geq r+1},F),\ell)
 +
\Delta^{r-1}P(S/(I,F)_{\geq r},\ell)}\\[2mm]
& \qquad\qquad= & \Delta^{r-1}(P(S/(I,F),\ell)-P(S/(I,F)_{\geq
 r},\ell))\\
& & -\Delta^{r-1}P(S/(I_{\geq r+1},F),\ell)+\Delta^{r-1}P(S/(I,F)_{\geq
r},\ell)\\[2mm]
& \qquad\qquad= & \Delta^{r-1}(\Delta_\tau
P(S/I,\ell)+P([0:F]_{S/I},\ell-\tau))\\
& & -\Delta^{r-1}\Delta_\tau P(S/I_{\geq r+1},\ell)\\[2mm]
& \qquad\qquad= & \Delta_\tau (\Delta^{r-1}(P(S/I,\ell)-P(S/I_{\geq
r+1},\ell)))\\
& & +\Delta^{r-1}P([0:F]_{S/I},\ell-\tau)\\[2mm]
& \qquad\qquad= & \tau\cdot\Delta^r(P(S/I,\ell)-P(S/I_{\geq r+1},\ell))
 +
\Delta^{r-1}P([0:F]_{S/I},\ell-\tau)\\[2mm]
& \qquad\qquad= & \tau\cdot
\mbox{arith-deg}_r(I)+\Delta^{r-1}P([0:F]_{S/I},\ell-\tau),
\end{eqnarray*}
by Lemma 2.5.  Hence the assertion is proved.
\vspace{5mm}

The following lemma is used in the proof of Theorem 2.1 and Lemma
 4.3.
\vspace{5mm}

\noindent{\bf Lemma 2.7:}{\em\
Let $I$ be a homogeneous ideal of $S$. Let $r$ be an integer. Let
 $F$ be a
homogeneous
polynomial of $S$ with degree $(F)\geq 1$ such that $F$ does not
 belong to
any associated
prime ideal $\p$ of $I$ with $\dim{\p}\geq r$. Then we
 have
$${(I_{\geq u},F)}_{\geq r} = (I,F)_{\geq r}$$
for all integers $u=-1,0,\cdots,r+1$.}
\vspace{5mm}

\noindent{\em Proof.}
We want to prove that
$${(I_{\geq u},F)}_{\p} = (I,F)_{\p}$$
for all prime ideals $\p$ with $\dim{\p}\geq r$. We may
 assume
that $F\in{\p}$.
If $\p$ does not contain any primary component $\q$ of
 $I$ with
$\dim({\q})
<u$, the proof is done. Now assume that there is a primary component
 $\q$ of $I$ with
$\dim({\q})<u$ such that ${\p}\rixrel{=}{\supset}{\q}$.
Since $r\leq\dim({\p})
\leq\dim({\q})<u$, we see that $u=r+1$ and ${\p}=\sqrt{\q}$.
Thus we
have that $\p$ is an associated prime ideal of $I$ with $\dim{\p}=r$ and
$F\in{\p}$, which contradicts the hypothesis.
\vspace{5mm}

\noindent{\em Proof of Theorem 2.1.}
First we prove the case $r\geq 1$. Applying (2.5) and (2.6), we
 have
\begin{eqnarray*}
\lefteqn{\mbox{arith-deg}_{r-1}(I,F)-\mbox{arith-deg}_{r-1}(I_{\geq
 r+1},F)}\\
& \quad= & \tau\cdot\mbox{arith-deg}_r(I) + \Delta^{r-1}P(S/(I_{\geq
r+1},F),\ell)\\
& & -\Delta^{r-1}P(S/(I,F)_{\geq r},\ell) +
\Delta^{r-1}P([0:F]_{S/I},\ell-\tau)\\
& & -\Delta^{r-1}(P(S/(I_{\geq
r+1},F),\ell)-P(S/(I_{\geq r+1},F)_{\geq r},\ell))\\
& \quad= & \tau\cdot\mbox{arith-deg}_r(I) + \Delta^{r-1}[P(S/(I_{\geq
r+1},F)_{\geq r},\ell)
-P(S/(I,F)_{\geq r},\ell)]\\
& & + \Delta^{r-1}P([0:F]_{S/I},\ell-\tau).
\end{eqnarray*}
By the assumption, $[0:F]_{S/I}$ has at most $(r-1)$-dimensional
 support,
which means
$\Delta^{r-1}P([0:F]_{S/I},\ell-\tau)\geq 0$. Further,
$\Delta^{r-1}P([0:F]_{S/I},\ell
-\tau)=0$ if and only if $F$ does not belong to any associated
 prime ideal
$\p$ with
$\dim{\p}=r-1$. On the other hand, $S/(I_{\geq r+1},F)_{\geq
r}=S/(I,F)_{\geq r}$
by Lemma 2.7. Therefore we have
$$\mbox{arith-deg}_{r-1}(I,F) - \mbox{arith-deg}_{r-1}(I_{\geq
r+1},F)\geq\tau\cdot
\mbox{arith-deg}_r(I)$$
and the equality holds if and only if $F$ does not belong to any
 associated
prime ideal
$\p$ with $\dim{\p}=r-1$.

Next we prove the case $r=0$. Now we see that
\begin{eqnarray*}
\mbox{arith-deg}_{-1}(I,F) & = & \mbox{length}_S(I,F)_{\geq 0}/(I,F)\\
& = & \sum_{\ell=0}^{N}(\dim_K{[S/(I,F)]}_\ell - \dim_K{[S/(I,F)_{\geq
0}]}_\ell)\\
& = & \sum_{\ell=0}^N(\dim_K{[S/I]}_\ell-\dim_K{[S/I]}_{\ell-\tau}\\
& & + \dim_K
{\left[[0:F]_{S/I}\right]}_{\ell-\tau}
    - \dim_K{[S/(I,F)_{\geq 0}]}_\ell)
\end{eqnarray*}
for large $N$. Similarly, we see that
\begin{eqnarray*}
\mbox{arith-deg}_{-1}(I_{\geq 1},F) & = &
\sum_{\ell=0}^N\left(\dim_K{[S/I_{\geq 1}]}_\ell
- \dim_K[S/I_{\geq 1}]_{\ell-\tau}\right.\\
& & \left. - \dim_K{[{S/(I_{\geq 1},F)}_{\geq 0}]}_\ell\right)
\end{eqnarray*}
for large $N$. Hence we have
\begin{eqnarray*}
\lefteqn{\mbox{arith-deg}_{-1}(I,F) - \mbox{arith-deg}_{-1}(I_{\geq
 1},F)}\\
& \qquad= & \sum_{\ell=0}^N(\dim_K{[I_{\geq 1}/I]}_\ell - \dim_K
 {[I_{\geq
1}/I]}_{\ell-\tau})\\
& & - \sum_{\ell=0}^N\dim_K {[{(I_{\geq 1},F)}_{\geq 0}/{(I,F)}_{\geq
 0}]}_\ell
+ \sum_{\ell=0}^N\dim_K{[{[0:F]}_{S/I}]}_{\ell-\tau}\\
& \qquad= & \sum_{\ell=N-\tau+1}^N\dim_K{[I_{\geq 1}/I]}_\ell
 -
\sum_{\ell=0}^N\dim_K
{[{(I_{\geq 1},F)}_{\geq 0}/{(I,F)}_{\geq 0}]}_\ell\\
& & + \sum_{\ell=0}^N\dim_K{[{[0:F]}_{S/I}]}_{\ell-\tau}
\end{eqnarray*}
for large $N$. By the assumption,
$\sum_{\ell=0}^N\dim_K{[{[0:F]}_{S/I}]}_{\ell-\tau}
=\mbox{length}_S({[0:F]}_{S/I})  \geq  0$ for large $N$, and ${[0:F]}_{S/I}
=0$
if and only if
$F$ is a non-zero-divisor of $S/I$. On the other hand, we see
$\dim_K{[I_{\geq 1}/I]}_\ell
=P(I_{\geq 1}/I,\ell)=\mbox{arith-deg}_0 (I)$ for large $\ell$.
 Further, we have
${(I_{\geq 0},F)}_{\p} = (I,F)_{\p}$ for all prime ideals
 $\p$ with
$\dim{\p}=0$ by Lemma 2.7. Hence we have
$$\mbox{arith-deg}_{-1}(I,F) -
\mbox{arith-deg}_{-1}(I_{\geq 1},F)\geq\tau\cdot
\mbox{arith-deg}_0 (I)$$
and the equality holds if and only if $F$ is a non-zero-divisor of
 $S/I$. This
completes the proof of Theorem 2.1.

\section{Castelnuovo-Mumford regularity}

Bayer and Mumford \cite{BM} give a bound for the arithmetic degree
in terms of the Castelnuovo-Mumford regularity. The aim of this section
is to describe improved bound on this degree.

Let $m = m(I)$ be the Castelnuovo-Mumford regularity (see, e.g.,
\cite{BM}, \cite{EG}, \cite{M}) for a homogeneous ideal $I$ of the
polynomial ring $S = K[x_0,\cdots,x_n]$. Then our main result is
the following theorem.

\vspace{5mm}

\noindent{\bf Theorem 3.1:}{\em\
Let $I$ be a homogeneous ideal of $S$.
Let $m = m(I)$ be the Castelnuovo-Mumford regularity of $I$.
Then we have, for any integer $r\geq 0$
$$\mbox{arith-deg}_r (I)\leq \Delta^r P(S/I,\ell)$$
for all integers $\ell\geq m-1$.}
\vspace{5mm}

We want to give two corollaries.  The first one shows
 that
(3.1) improves
the bound given in \cite{BM}, Proposition 3.6.
\vspace{5mm}

\noindent{\bf Corollary 3.2:}{\em\
For all $r \geq 0$, we have
\begin{eqnarray*}
\mbox{arith-deg}_r (I)
& \leq & \Delta^r P(S/I,m-1)\\
& \leq & \left(\begin{array}{c}m+n-r-1\\n-r\end{array}\right)\\
& \leq & m^{n-r}
\end{eqnarray*}}

\noindent{\bf Corollary 3.3:}{\em\
Let $t=\mbox{depth } S/I$. Then we have, for an integer
 $r\geq 0$,
$$\mbox{arith-deg}_r (I)\leq\Delta^r H(S/I,\ell)$$
for all $\ell\geq m+r-t-1$ if $r-t$ is even, and for all $\ell\geq
 m+r-t$
if $r-t$ is
odd.}
\vspace{5mm}

Before embarking on the proof of Theorem 3.1 and the corollaries
 we need
two
lemmas. The first one follows from \cite{S} Nr.79 (see also \cite{SV},
Proof of Lemma
I.4.3).

\vspace{5mm}

\noindent{\bf Lemma 3.4:}{\em\
 Let $I$ be a homogeneous ideal of $S$ and
 $t=\mbox{depth } S/I$.
 Then we have
\renewcommand{\labelenumi}{(\theenumi)}
\renewcommand{\theenumi}{\alph{enumi}}
\begin{enumerate}
\item $P(S/I,\ell) = H(S/I,\ell) - \sum_{i=0}^d (-1)^i\dim_K
 {[H_{\m}^i(S/I)]}_\ell$
for all integers $\ell$.
\item $P(S/I,\ell) = H(S/I,\ell)$ for all $\ell\geq m-t$.

\end{enumerate}
}
\vspace{5mm}

\noindent{\bf Lemma 3.5:}{\em\
Let $I$ be a homogeneous ideal of $S$.
Then we have
$$\Delta^r P(I,\ell)\geq 0$$
\vspace{-10mm}

\noindent for all $\ell\geq m-1$ and $r\geq 0$.
}
\vspace{5mm}

\noindent{\em Proof.}
For a generic hyperplane $H$ given by $h=0$, we can take an exact
 sequence
$$0\rightarrow I(-1)\stackrel{h}{\rightarrow}
I\rightarrow I_H\rightarrow 0 ,$$
where $I_H=(I,h)/h$.
  From the exact sequence, we have
$$\Delta P(I,\ell) = P(I_H,\ell)$$
for all $\ell$ and $I_H$ is $m$-regular. Repeating this step, we
 see that
$$\Delta^r P(I,\ell) = P(I_{H_1\cap\cdots\cap H_r},\ell)$$
for all $\ell$ and for generic hyperplanes $H_1,\cdots,H_r$ defined
 by
$h_1=0,\cdots,h_r=0$, resp., where
$$I_{H_1\cap\cdots\cap H_r} = (I,h_1,\cdots,h_r)/(h_1,\cdots,h_r),$$
 and that
$I_{H_1\cap\cdots\cap H_r}$ is $m$-regular. So ${(I_{H_1\cap\cdots\cap
H_r})}_{\geq 0}$
is also $m$-regular and \linebreak $\mbox{depth}_S S/{(I_{H_1\cap\cdots\cap
H_r})}_{\geq 0}\geq 1$.
Therefore we have
\begin{eqnarray*}
\Delta^r P(I,\ell) & = & P(I_{H_1\cap\cdots\cap H_r},\ell)\\
& = & P(S,\ell) - P(S/I_{H_1\cap\cdots\cap H_r},\ell)\\
& = & P(S,\ell) - P(S/{(I_{H_1\cap\cdots\cap H_r})}_{\geq 0},\ell)\\
& = & H(S,\ell) - H(S/{(I_{H_1\cap\cdots\cap H_r})}_{\geq 0},\ell)\\
& = & H({(I_{H_1\cap\cdots\cap H_r})}_{\geq 0},\ell)\geq 0
\end{eqnarray*}
for $\ell\geq m-1$, by (3.4), (b).
\vspace{5mm}

\noindent{\em Proof of Theorem 3.1.}
Without the loss of generality, we may assume\linebreak that $I$
 is a saturated
ideal. First we
prove the case $r=0$. By Lemma 2.5, we have
$$(\ast)\qquad\qquad\mbox{arith-deg}_0(I)=P(S/I,\ell) - P(S/I_{\geq
 1},\ell).$$
Now we want to show that $I_{\geq 1}$ is $m$-regular. From the short
 exact
sequence
$$0\rightarrow I_{\geq 1}/I\rightarrow S/I\rightarrow S/I_{\geq
1}\rightarrow 0$$
and the fact that $I_{\geq 1}/I$ has at most 1-dimensional support and
by Grothendieck's vanishing theorem $H_{\m}^i(I_{\ge 1}/I)  =  0$
for $i  \ge  2$,
 we have
$$0\rightarrow H_{\m}^1 (I_{\geq 1}/I)\rightarrow
H_{\m}^1(S/I)\rightarrow H^1_{\m}
(S/I_{\geq 1})\rightarrow 0$$
and $H_{\m}^i(S/I)\cong H_{\m}^i (S/I_{\geq 1})$ for $i\geq
 2$.
Thus we have $I_{\geq 1}$
is $m$-regular. Hence
$$P(S/I_{\geq 1},\ell) = H(S/I_{\geq 1},\ell)\geq 0$$
for all $\ell\geq m-1$, by (3.4), (b). Therefore we have from $(\ast)$
$$\mbox{arith-deg}_0 (I)\leq P(S/I,\ell)$$
for all $\ell\geq m-1$.

Now let us assume $r>0$. By Corollary 2.4 we see
$$\mbox{arith-deg}_r (I) = \mbox{arith-deg}_0 (I,h_1,\cdots,h_r)$$
for generic hyperplanes $h_1,\cdots,h_r$. Thus we have
$$\mbox{arith-deg}_r (I)\leq P(S/(I,h_1,\cdots,h_r),\ell)$$
for all $\ell\geq m-1$. On the other hand, we see
\begin{eqnarray*}
P(S/(I,h_1,\cdots,h_r),\ell) & = & \Delta P(S/(I,h_1,\cdots,h_{r-1}),\ell)
\\
& & \qquad \vdots\\
& = & \Delta^r P(S/I,\ell)
\end{eqnarray*}
for all $\ell$. Hence the assertion is proved.
\vspace{5mm}

\noindent{\em Proof of Corollary 3.2.}
By Lemma 3.5, we have
$$\Delta^r P(S/I,\ell)\leq\Delta^r P(S,\ell)$$
for all $\ell\geq m-1$. On the other hand,
$\Delta^r P(S,\ell) =
\left(\begin{array}{c}
n+\ell-r\\n-r\end{array}\right)$. Hence the assertion
follows from Theorem 3.1.
\vspace{5mm}

\noindent{\em Proof of Corollary 3.3.}
By Lemma 3.4, we see that
$$\Delta^r P(S/I,\ell) = \Delta^r H(S/I,\ell)$$
for all $\ell\geq m+r-t$, and that
$$\Delta^r P(S/I,m+r-t-1)=\Delta^r H(S/I,m+r-t-1)-(-1)^{r-t}\dim_K
{[H_{\m}^t(S/I)]}_{m-t-1}$$
because
$$P(S/I,m-t-1)  =  H(S/I,m-t-1)  -  (-1)^t[H_{\m}^t(S/I)]_{m-t-1}.$$
Hence the assertion follows from Theorem 3.1.

\section{Bezout-type results}

The aim of this section is to state properties of arithmetic degree
 under
iterated
hyperplane sections, and Bezout-type results. Our Theorem 4.1 describes
a Bezout's theorem in terms of the arithmetic degree.
\vspace{5mm}

\noindent{\bf Theorem 4.1:}{\em\
Let $I$ be a homogeneous ideal of $S := K[x_0,x_1,\cdots,x_n]$.
 Let $r\geq
0$ and
$s\geq 1$ be integers with $r+1\geq s$. Let $F_1,\cdots,F_s$ be homogeneous
polynomials
of $S$ such that $F_i$ does not belong to any associated prime ideal
${\p}$ of
$(I,F_1,\cdots,F_{i-1})$ with $\dim {\p}\geq r-i+1$, for all
$i=1,\cdots,s$. Then we have
\renewcommand{\labelenumi}{(\theenumi)}
\renewcommand{\theenumi}{\roman{enumi}}
\begin{enumerate}
\item $\mbox{arith-deg}_{r-s} (I,F_1,\cdots,F_s)\geq\left[\prod_{i=1}^s
\mbox{degree} (F_i)\right]\cdot\mbox{arith-deg}_r (I)$;
\item We have equality in (i) if and only if
$({(I,F_1,\cdots,F_{i-1})}_{\geq r-i+2},F_i)$
has no $(r-i)$-dimensional primes and $F_i$ does not belong to any
associated prime ideal
${\p}$ of $(I,F_1,\cdots,F_{i-1})$ with $\dim {\p}=r-i$,
 for all
$i=1,\cdots,s$;
\item Assume that there is an integer $t$ with $-1\leq t\leq r+1$
 such that
$F_i$ does
not belong to any associated prime ideal ${\p}$ of $(I_{\geq
t},F_1,\cdots,F_{i-1})$ with
$\dim {\p}\geq r-i+1$, for all $i=1,\cdots,s$, and $(I_{\geq
t},F_1,\cdots,F_s)$ has no
$(r-s)$-dimensional associated prime ideals. Then we have equality
 in (i)
if and only if
$F_i$ does not belong to any associated prime ideal ${\p}$ of
$(I,F_1,\cdots,F_{i-1})$
with $\dim {\p}=r-i$, for all $i=1,\cdots,s$.
\end{enumerate}}
\vspace{5mm}

\noindent{\em Proof.}
(i) and (ii) follow from (2.2).
In order to prove (iii) we need Lemma 4.2 and Lemma 4.3 below. First
 we
replace the ideal
$I$ of (4.2) by the ideal $I_{\geq t}$ of (iii). Then Lemma 4.3 shows
 that
we can apply
(ii) of (4.1). This provides our result (iii).
\vspace{5mm}

We note that special cases of (4.1) describe
 generalizations of
classical results in the degree theory (see, e.g., \cite{FV}, \cite{V}).
\vspace{5mm}

We prove the two lemmas.
\vspace{5mm}

\noindent{\bf Lemma 4.2:}{\em\
Let $I$ be a homogeneous ideal of $S$. Let $r$ and $s$ be integers
 with
$1 \leq s \leq r+1$. Let $F_1, \cdots, F_s$ be homogeneous polynomials
 of $S$
with\linebreak
degree $(F_i) \geq 1$, $i=1,\cdots,s$, such that $F_i$ does not
 belong to
any associated
prime ideal $\p$ of $(I,F_1,\cdots,F_{i-1})$ with $\dim{\p}\geq
r-i+1$, for all
$i=1,\cdots,s$. If the ideal $(I,F_1,\cdots,F_s)$ has no
$(r-s)$-dimensional associated
prime ideals, then the ideal $({(I,F_1,\cdots,F_{i-1})}_{\geq
r-i+2},F_i)$ has no
$(r-i)$-dimensional associated prime ideals, for all $i=1,\cdots,s$.}
\vspace{5mm}

\noindent{\em Proof.}
It is easy to see that the ideal $(I,F_1,\cdots,
F_i)$
has no
$(r-i)$-dimensional associated prime ideals. By Theorem 2.1, we have
$$\mbox{arith-deg}_{r-i} {(I,F_1,\cdots,F_i)}_{\geq r-i} -
\mbox{arith-deg}_{r-i}
({(I,F_1,\cdots,F_{i-1})}_{\geq r-i+2},F_i)\geq 0.$$
This shows that $\mbox{arith-deg}_{r-i}({(I,F_1,\cdots,F_{i-1})}_{\geq
r-i+2},F_i)=0$.
Thus the assertion is proved.
\vspace{5mm}

\noindent{\bf Lemma 4.3:}{\em\
Let $I$ be a homogeneous ideal of $S$. Let $r$ and $s$ be integers
 with
$1\leq s\leq r+1$.
Let $t$ be an integer with $-1\leq t\leq r+1$. Let $F_1,\cdots,F_s$
 be
homogeneous
polynomials of $S$ with degree $(F_i)\geq 1$, $i=1,\cdots,s$, such
 that
$F_i$ does not
belong to any associated prime ideal $\p$ of $(I,F_1,\cdots,F_{i-1})$
 with
$\dim{\p}\geq r-i+1$, for all $i=1,\cdots,s$. Then we have
$${(I_{\geq t},F_1,\cdots,F_{i-1})}_{\geq r-i+2}=
{(I,F_1,\cdots,F_{i-1})}_{\geq r-i+2}$$
for $i=1,\cdots,s$.}
\vspace{5mm}

\noindent{\em Proof.}
By Lemma 2.7, we have
\begin{eqnarray*}
{(I,F_1,\cdots,F_{i-1})}_{\geq r-i+2} & = &
{({(I,F_1,\cdots,F_{i-2})}_{\geq r-i+3},
F_{i-1})}_{\geq r-i+2}\\
& = & \cdots=\ {(\cdots{({(I,F_1)}_{\geq r},F_2)}_{\geq
r-1},\cdots,F_{i-1})}_{\geq r-i+2}\\
& = & {(\cdots{({(I_{\geq t},F_1)}_{\geq r},F_2)}_{\geq
r-1},\cdots,F_{i-1})}_{\geq r-i+2}.
\end{eqnarray*}
On the other hand, we see
$${(I_{\geq t},F_1,F_2,\cdots,F_{i-1})}_{\geq r-i+2}\subset
{(\cdots{({(I_{\geq t},F_1)}_{\geq r},F_2)}_{\geq
r-1},\cdots,F_{i-1})}_{\geq r-i+2}$$
and
$${(I_{\geq t},F_1,F_2,\cdots,F_{i-1})}_{\geq r-i+2}\supset
{(I,F_1,\cdots,F_{i-2})}_{\geq r-i+2}.$$
Therefore we have ${(I_{\geq t},F_1,\cdots,F_{i-1})}_{\geq r-i+2}
 =
{(I,F_1,\cdots,F_{i-1})}_{\geq r-i+2}$ for all\linebreak $i=1,\cdots,s$.
\vspace{5mm}

\section{Some examples}

The first example sheds some light on Theorem 2.1 and Corollary 2.2
 in case that
$F$ has degree one and is a non-zero-divisor on $S/I$. It shows that
 we have no
equality in Corollary 2.2 even under these assumptions.
\vspace{5mm}

\noindent{\em Example 1:}
Let $S=K[x_0,x_1,x_2,x_3,y_1,y_2,\cdots,y_r]$ be a polynomial ring,
 where
$r$ is a
non-negative integer. Take
${\q}=(x_0x_3-x_1x_2,x_0^2,x_1^2,x_0x_1)\subset S$, which
is a primary ideal belonging to $(x_0,x_1)$ (cf. \cite{SV}, Claim
 1 on page
182).
We set $I={\q}\cap (x_0^2,x_1,x_2)$ and
$F(x_0,x_1,x_2,x_3)=x_3+G(x_0,x_1,x_2)$, where $G(x_0,x_1,x_2)$
is a linear form. Then we have by (2.1)
$$\mbox{arith-deg}_{r-1}(I,F) - \mbox{arith-deg}_{r-1}(I_{\geq r+1},F)
=\mbox{arith-deg}_r (I).$$
Now we will show that
$$\mbox{arith-deg}_{r-1}(I,F)>\mbox{arith-deg}_r (I).$$
For simplicity we assume that $G=0$, that is, $F=x_3$. Clearly,
$\mbox{arith-deg}_r (I)\linebreak=1$.
On the other hand,
\begin{eqnarray*}
(I,x_3) & = & (x_0^2,x_1^2,x_0x_1,x_0x_2x_3-x_1x_2^2) + (x_3) \\
& = & (x_0^2,x_1,x_3)\cap (x_0^2,x_1^2,x_2^2,x_0x_1,x_3).
\end{eqnarray*}
Hence $\mbox{arith-deg}_{r-1}(I,x_3)=2$. Also, we have
\begin{eqnarray*}
(I_{\geq r+1},x_3) & = & ({\q},x_3)\\
& = & (x_0^2,x_1^2,x_0x_1,x_1x_2,x_3)\\
& = & (x_0^2,x_1,x_3)\cap (x_0^2,x_1^2,x_2,x_3,x_0x_1).
\end{eqnarray*}
Hence $\mbox{arith-deg}_{r-1}(I_{\geq r+1},x_3)=1$.

We note that
$$\mbox{arith-deg}_{r-1}(I,x_3)>\mbox{arith-deg}_r (I)$$
even in the case that $x_3$ is a non-zero-divisor on $S/I$.
\vspace{5mm}

The second example shows that the bound of Theorem 3.1 is sharp and
 improves the
result of \cite{BM}, Proposition 3.6 (see Corollary 3.2).
\vspace{5mm}

\noindent{\em Example 2:}
Take
$I=(x_0^2x_1,x_0x_2^2,x_1^2,x_2,x_2^3)
=(x_0^2,x_2)\cap(x_1,x_2^2)\cap(x_0^2,
x_1^2,
x_0x_2^2,x_2^3)\linebreak\subset S := K[x_0,x_1,x_2]$.
We get $m=5$,
$P(S/I,\ell)=4$ for all $\ell\geq 4$. We consider the case $r=0$
 in (3.1)
and (3.2).
Then we have
$$4 = \deg I = \mbox{arith-deg}_0 (I)\leq P(S/I,m-1)=4
<\left(\begin{array}{c}
5+2-0-1\\2-0\end{array}\right)=15.$$

\end{document}